\documentclass[a4paper,12pt]{article}
\usepackage{amsmath, amsfonts, amssymb,amsthm, mathrsfs}
\usepackage[english]{babel}

\begin{document}

\begin{center} \bf Arens algebras of measurable operators for Maharam
traces \\
A.A.Alimov \\
e-mail: alimovakrom63@yandex.ru
\end{center}

\date{}
\begin{abstract}
We study order and topological properties of
the non-commutative Arens algebra   associated with arbitrary
Maharam trace.
\end{abstract}



\section{Introduction}

The integration theory for traces and weights, given on von
Neumann algebras, is one of the central objects for numerous
investigations connected with operator algebras and their
applications. Development of the non-commutative theory begun with
the work by I.~Segal \cite{S11}, where the author introduced the
$\ast$-algebra $S(M)$ of all measurable operators affiliated with
a von Neumann algebra $M$ being a non-commutative analog of the
$\ast$-algebra $L^{0}(\Omega, \Sigma, \mu)$ of all measurable
complex functions given on a measurable space $(\Omega, \Sigma,
\mu).$ This algebra $S(M)$ has became the base for construction of
the general theory of non-commutative $L^{p}$ spaces associated
with  a faithful normal semi-finite trace $\mu$ given on $M$ (see,
for example, \cite{Y11, FK1, PX1}. The Banach spaces $L^{p}(M,
\tau), p\geq 1,$ are ideal linear subspaces in the $\ast$-algebra
$S(M, \tau)$ of all $\tau$-measurable operators affiliated with
$M.$ The $\ast$-algebra $S(M, \tau)$, introduced in \cite{N11}, is
itself a $\ast$-subalgebra in $S(M)$ and it coincides with the
completion of $M$ with respect to the topology of convergence in
measure generated by the trace $\tau.$ It should be noted, both
$^\ast$-algebras $S(M)$ and $S(M, \tau)$ are used for meaningful
examples of $EW^{\ast}$-algebras of unbounded operators that are
important in the general theory of algebras of unbounded operators
(see, for example \cite{AI1}). Interesting examples of
$EW^{\ast}$-algebras are also the Arens $^\ast$-algebras
$L^{\omega}(M, \tau) = \bigcap\limits_{p\geq 1} L^{p}(M, \tau).$
The Arens algebras $L^{\omega}(M, \tau)$ were considered at first
in \cite{A11} for the case of $M = L^{\infty}(0,1).$
Non-commutative Arens algebras were introduced by A.~Inoue in
\cite{I11}, and then properties of such algebras were studied in
\cite{A12,Z11}.

Due to presence of center-valued traces on finite von Neumann
algebras, it is natural to extend the integration theory for
traces with values in a complex order-complete lattices
$F_{\mathbb{C}} = F\oplus iF.$

If the original von Neumann algebra is commutative, then
constructing of $F_{\mathbb{C}}$-valued integration for it is a
component part of investigations of properties of order-continuous
positive mappings of vector lattices. The theory of such mappings
was described in details in \cite{K11}, chapter 4. Operators,
having the Maharam property, are important among these mappings.
$L^{p}$-spaces associated with such operators are significant
examples of Banach-Kantorovich vector lattices.

For non-commutative von Neumann algebras $M,$ properties of spaces
$(L^{p}(M, \Phi),\parallel\cdot\parallel_{p})$, constructed by a
$F_{\mathbb{C}}$-valued trace $\Phi$, were considered in 
\cite{GCH} and \cite{CHK} in the case $F_{\mathbb{C}}$ is a von
Neumann subalgebra in the center of the algebra, and the trace
$\Phi$ has the following modularity property: $\Phi (zx) = z\Phi
(x)$ for all $z\in F_{\mathbb{C}}, x\in M.$

The modularity property implies immediately that if $0\leq
f\leq\Phi(x),$ $f\in F$, $ x\in M,$ then there exists $ y\in M,$ $
0\leq y\leq x$ such that $\Phi(y)=f.$ It means that the trace
$\Phi$ possesses the Maharam property (compare \cite{K11}, 3.4.1).

Faithful normal traces $\Phi$ on a von Neumann algebra $M$ with
values in arbitrary complex order-complete vector lattice were
considered in \cite{CHZ}, where, in particular, full description
of such traces is given for the case when $\Phi$ is the Maharam
trace. In \cite{CHZ}, the non-commutative $L^{1}$-space $L^{1}(M,
\Phi)\subset S(M)$ associated with the Maharam trace $\Phi$ was
constructed with the help of the in topology of convergence the
measure, and it was established that $L^{1}(M, \Phi)$ is a
Banach-Kantorovich space. Later in \cite{CZ1}, non-commutative
$L^{p}$-spaces $L^{p}(M, \Phi)$ were defined for all $p>1.$ The
problem of construction and description of properties of the Arens
algebras $L^{\omega}(M, \Phi),$ associated with the Maharam trace
$\Phi$ has arisen naturally. Such problem is solved in \cite{A13}
for the case when $\Phi$ takes values in $S(A),$ where $A$ is a
von Neumann subalgebra in the center $Z(M)$ of an algebra $M.$

In the present paper, we study order and topological properties of
the Arens algebra  $L^{\omega}(M, \Phi)$ associated with arbitrary
Maharam trace $\Phi.$

Necessary and sufficient conditions are determined that provide
local convexity of the topology $\tau_{\omega}(M,\Phi)$ in
$L^{\omega}(M,\Phi)$ generated by the system of norms
$\{\parallel\cdot\parallel_{p}\}_{p\geq 1}$. A criterion for
coincidence of the topology $\tau_{\omega}(M,\Phi)$ with the
$(o)$-topology in the ordered linear space $L_{h}^{\omega}(M,
\Phi)= \{x\in L^{\omega}(M,\Phi):x=x^{\ast}\}$ is established

We use terminology and notations of the theory of von Neumann
algebras from \cite{SZ1, T11}, the theory of measurable operators
from \cite{S11, MCH} and the theory of vector lattices from
\cite{K11}.

\section{Preliminaries}

Let $H$ be a Hilbert space over the field of complex numbers
$\mathbb{C},$ $B(H)$ be a $\ast$-algebra of all bounded linear
operators acting in $H,$ $\mathbf{1}$ be an identical operator in
$H$, and let $M$ be a von Neumann subalgebra in $B(H)$. Denote by
$P(M) = \{ p\in M: \ p^{2}=p=p^{\ast} \}$ the lattice of all
projects from $M,$ and by $P_{fin}(M)$ the sublattice of all
finite projections from $M.$

A closed linear operator $x,$ affiliated with the von Neumann
algebra $M,$ having a dense domain  $D(x)\subset H,$ is called
\textit{measurable} with respect to $M$ if there exists a sequence
$\{p_{n}\}_{n=1}^{\infty}\subset P(M)$ such that $p_{n}\uparrow
\mathbf{1},$ $p_{n}(H)\subset D(x)$ and $p_{n}^{\perp} =
\mathbf{1} - p_{n}\in P_{fin}(M)$ for each $n=1, 2, \dots$

The set $S(M)$ of all measurable operators with respect to $M$ is
a $\ast$-algebra with the unit $\mathbf{1}$ over $\mathbb{C}$ with
respect to the natural involution, multiplication on a scalar, and
operations of the strong addition and strong multiplication
obtained by closure of usual operations \cite{S11}. It is clear,
$M$ is a $\ast$-subalgebra in $S(M).$

If $x\in S(M),$ and $x = u\mid x\mid$ is the polar decomposition
of the operator $x$ where $\mid x\mid =
(x^{\ast}x)^{\frac{1}{2}},$ $u$ is the corresponding partial
isometry from $B(H)$ for which $u^{\ast}u$ is the right support
for $x,$ then $u\in M$ and $\mid x\mid \in S(M).$ The spectral
family of projectors $\{E_{\lambda}(x)\}_{\lambda\in\mathbb{R}}$
of the self-adjoint operator $x\in S(M)$ is always contained in
$P(M).$

For any subset $E\subset S(M),$ we denote by $E_{h}$ ($E_{+}$,
respectively) the set of all self-adjoint (positive, respectively)
operators from $E$.

Let $M$ be a commutative von Neumann algebra. In this case, there
exists a faithful normal semi-finite trace $\tau$ on $M,$ and $M$
is $^*$-isomorphic to $W^*$-algebra $L^{\infty}(\Omega, \Sigma,
\mu)$ of all essentially bounded complex measurable functions
given on a measurable space $(\Omega, \Sigma, \mu)$ with the
locally finite measure $\mu$ having the direct sum property
(almost everywhere  equal functions are identified). Moreover,
$\mu (A) = \tau(\tilde{\chi}_A),$ $A\in\Sigma$ where
$\tilde{\chi}_A$ is the equivalence class containing the function
${\chi}_A$ (recall that ${\chi}_{A}(\omega) = 1$ for $\omega\in
A,$ and ${\chi}_{A} (\omega) = 0$ if $\omega\notin A ).$ In
addition, the $^*$-algebra $S(M)$ is identified with the
$^*$-algebra $L^{0}(\Omega, \Sigma, \mu)$ of all measurable
complex functions given on $(\Omega, \Sigma, \mu)$ (almost
everywhere equal functions are identified) \cite{S11}. Let us
consider in $L^{0}(\Omega, \Sigma, \mu)$ the topology $t(M)$ of
the convergence locally in the measure, i.e. the Hausdorff
topology endowing $L^{0}(\Omega, \Sigma, \mu)$ with the structure
of a complete topological $^*$-algebra, the base of zero
neighborhoods of which is formed by the sets in the form of
$$
W(B, \varepsilon, \delta) = \left\{f\in L^{0}(\Omega, \Sigma,
\mu): \ \mbox{ there exists a set } \ E\in\Sigma\ \mbox{ such that
} \right.
$$
\begin{equation}
\left. E\subseteq B, \mu (B\setminus E)\leq\delta, f\chi_{E}\in
L^{\infty}(\Omega, \Sigma, \mu), \|f\chi_{E}\|_{L^{\infty}(\Omega,
\Sigma, \mu)}\leq\varepsilon\right\},
\end{equation}
where $\varepsilon, \delta > 0,$ $B\in\Sigma,$ $\mu (B)<\infty,$
$\|\cdot\|_{L^{\infty}(\Omega, \Sigma, \mu)}$ is the $C^{*}$-norm
in $L^{\infty}(\Omega, \Sigma, \mu).$ The sets
$W(B,\varepsilon,\delta)$ have the following ideality property: if
$g\in L^{0}(\Omega,\Sigma,\mu),$ \ $f\in W(B,\varepsilon,\delta),$
and $|g|\leq|f|$, then $g\in W(B,\varepsilon, \delta).$

Convergence of the net $f_{\alpha}$ to $f$ in the topology $t(M)$
(notation: $f_{\alpha}\xrightarrow {t(M)}f$) means that
$f_{\alpha}\chi_{B}\rightarrow f\chi_{B}$ by the measure $\mu$ for
any $B\in\Sigma$ with $\mu(B)<\infty.$ Evidently, the topology
$t(M)$ is not changed if a trace $\tau$ is replaced by another
faithful normal semi-finite trace on $M.$ Therefore the topology
is uniquely defined by the von Neumann algebra $M$ itself. It  is
clear that the topology $t(M)$ is metrizable if and only if the
algebra $M$ is $\sigma$-finite, i.e. any set of nonzero mutually
orthogonal projections at most countable.

Now let $M$ be an arbitrary finite von Neumann algebra, $Z(M)$ be
the center in $M,$ and $\Phi_{M}: M\rightarrow Z(M)$ be a
center-valued trace on $M$ (\cite{SZ1}, 7.11). Identify the center
$Z(M)$ with the $\ast$-algebra $L^{\infty}(\Omega, \Sigma, \mu)$
and $S(Z(M))$ with the $\ast$-algebra $L^{0}(\Omega, \Sigma,
\mu).$ For arbitrary numbers $\varepsilon, \delta >0$ and
arbitrary set $B\in\Sigma$ with the measure $\mu (B)<\infty,$ set:
\begin{center} $V(B, \varepsilon, \delta) = \{x\in S(M):$ \
there exist $p\in P(M), z\in P(Z(M))$ such that $xp\in M,
\|xp\|_{M}\leq\varepsilon, z^{\perp}\in W(B, \varepsilon,
\delta),$ and $\Phi_{M}(zp^{\perp})\leq\varepsilon z \}$
\end{center} where $\|\cdot\|_{M}$ is the $C^{\ast}$-norm in $M.$
In \cite{MCH},$\S$ 3.5, it is shown that the system of sets
\begin{equation}
\{ \{x+V(B, \varepsilon, \delta)\}: \ x\in S(M), \
\varepsilon, \delta >0, \ B\in\Sigma, \ \mu(B)<\infty\}
\end{equation}
defines in $S(M)$ a Hausdorff vector topology $t(M),$ in which the
sets $(2)$ form the base of neighborhoods for the operator $x\in
S(M).$ In addition, $(S(M), t(M))$ is a complete topological
$^*$-algebra. The topology $t(M)$ is called the \textit{topology
of  convergence locally in a measure} [1]. It is clear, the
topology $t(M)$ induces in $S(Z(M))$ the topology $t(Z(M)).$
Moreover, if $Z(M)$ is a $\sigma$-finite algebra, then $t(M)$ is
metrizable.

The following criterion for convergence of nets in the topology $t(M)$ follows from
\cite{MCH}, $\S$ 3.5.\\

\noindent\textbf{Proposition 2.1.} {\it The net $\{x_{\alpha}\}_{\alpha\in
A}\subset S(M)$ converges to zero in the topology $t(M)$ if and only if
$\Phi_{M}(E_{\lambda}^{\perp}(\mid x_{\alpha}\mid))\xrightarrow{t(M)}0$
for any} $\lambda >0.$ \\

Let $F$ be an order complete vector lattice, $F_{\mathbb{C}} =
F\oplus iF$ be the complexification of $F,$ where $i$ is the
imaginary unit. As usual, for an element $z = \alpha + i\beta\in
F_{\mathbb{C}},$ $\alpha, \beta\in F,$ the adjoint element is
defined as $\bar{z} = \alpha - i\beta,$ and the module $\mid
z\mid$ is defined as $\mid z\mid : = \sup\{ Re(e^{i\theta}z):
0\leq\theta < 2\pi\}$ (\cite{K11}, 1.3.13).

A linear mapping $\Phi$ from a von Neumann algebra
$M$ into $F_{\mathbb{C}} $ is said to be an $F_{\mathbb{C}}$-valued trace if
$\Phi(x^{*}x) = \Phi(xx^{*})\geq 0$ for all $x\in M.$ It is clear, $\Phi (M_{h})\subset F,$
$\Phi (M_{+})\subset F_{+} = \{ a\in F, \ a\geq 0\}. $

The trace $\Phi$ is said to be {\it faithful} if the equality
$\Phi(x^{*}x) = 0$ implies $x=0.$ As well as for numerical traces
(see, for example,(\cite{T11}, chapter $V, \S 2$), it is
established that if there is a faithful $F_{\mathbb{C}}$-valued
trace on a von Neumann algebra $M,$ then this algebra is finite.
Let us list
some necessary properties of faithful traces $\Phi : M\rightarrow F_{\mathbb{C}}.$ \\

\noindent \textbf{Proposition 2.2} (\cite{ZCH}). {\it For any $x,
y, a, b\in M$ the following relations hold:

$1.$ $\Phi(x^*) = \overline{\Phi(x)};$

$2.$ $\Phi (xy) = \Phi (yx);$

$3.$ $\Phi(\mid x^*\mid) =  \Phi(\mid x\mid);$

$4.$ $\mid \Phi(axb)\mid\leq\parallel a\parallel_{M}\parallel
b\parallel_{M}\Phi(\mid x\mid);$

$5.$ If $x_{n},$ $x\in M,$ and $\parallel x_{n}- x\parallel_{M}
\rightarrow 0,$ then $\mid \Phi (x_{n}) - \Phi (x)\mid$ converges
in $F$ to zero with the regulator} $\Phi (\mathbf{1});$

$6.$ $\Phi (\mid x+y\mid)\leq \Phi (\mid x\mid) + \Phi (\mid
y\mid).$\\

We say that a trace $\Phi : M\rightarrow F_{\mathbb{C}} $
possesses the Maharam property if for any $x\in M_{+},$ $0\leq
f\leq\Phi (x),$ $f\in F$ there exists $y\in M_{+}$ such that
$y\leq x$ and $\Phi (y) =f.$

A trace $\Phi$ is said to be {\it normal} if $x_{\alpha},$ $x\in
M_{h},$ $x_{\alpha} \uparrow x$ imply $\Phi (x_{\alpha})\uparrow
\Phi (x).$ A faithful normal $F_{\mathbb{C}}$-valued trace $\Phi$
possessing the Maharam property is said to be {\it the Maharam
trace} (\cite{ZCH, CHZ}).

Let $\mathbf{1}_F$ be a weak unit in $F.$ Denote by $B(F)$ the
complete Boolean algebra of unit elements in $F$ with respect to
$\mathbf{1}_F.$ Let $Q$ be the Stone compact for $B(F),$ and let
$C_{\infty}(Q)$ be an extended order complete vector lattice of
all continuous functions $f: Q\rightarrow [-\infty; +\infty]$
taking the values $\pm\infty$  on nowhere dense sets from $Q.$
Identify $F$ with the fundament in the lattice $C_{\infty}(Q)$
consisting of the algebra $C(Q)$ of all continuous functions on
$Q,$ in addition $\mathbf{1}_F$ is identified with the function
which is identically equal to the unit (\cite{K11}, 1.4.4).

The following theorem from \cite{CHZ} gives description of Maharam
traces.\\

\noindent \textbf{Theorem 2.3} (\cite{CHZ}). {\it Let $\Phi$ be a
$F_{\mathbb{C}}$-valued trace on a von Neumann algebra $M.$ Then
there exist a von Neumann subalgebra $A$ in the center of $Z(M),$
an $\ast $-isomorphism  $\psi$ from $A$ onto the $\ast$-algebra
$C(Q)_{\mathbb{C}} = C(Q)\oplus iC(Q),$ a positive normal linear
operator $\mathscr{E}$ from $Z(M)$ onto $A$ with
$\mathscr{E}(\mathbf{1}) = (\mathbf{1}),$ $\mathscr{E}^{2} =
\mathscr{E},$ such that

$1.$ $\Phi (x) = \Phi (\mathbf{1})\cdot \psi (\mathscr{E}(\Phi
_{M}(x)))$ for all $x\in M;$

$2.$ $ \Phi (zy) = \Phi (z\mathscr{E}(y))$ for all $z,y\in Z(M);$

$3.$ $\Phi (zy) = \psi (z)\Phi (y)$ for all} $z\in A,$ $y\in M.$
\\

Theorem 2.3 implies that the $\ast$-algebra $B =
C(Q)_{\mathbb{C}}$ is a commutative von Neumann algebra and the
$\ast$-algebra $(C_{\infty}(Q))_{\mathbb{C}}$ is identified with
the $\ast$-algebra $S(B).$ We denote the unit of the algebra $B$
by $\mathbf{1}_{B}$ (it coincides with the weak unit
$\mathbf{1}_{F}$).

Let $\Phi$ be  $S(B)$-valued Maharam trace on the von Neumann
algebra  $M.$ Recall definition of the space $L^{1}(M, \Phi)$
[11]. We say that a net $\{x_{\alpha}\}\subset S(M)$ converges to
$x\in S(M)$ by a trace $\Phi$ (notation: $x_{\alpha}
\xrightarrow{\Phi} x)$ if $$\Phi (E_{\lambda}^{\perp}(\mid
x_{\alpha} -x\mid))\xrightarrow{t(B)} 0$$ for $\lambda >0.$ It was
shown in [11] that $x_{\alpha}\xrightarrow{\Phi} x \
\Leftrightarrow x_{\alpha}\xrightarrow{t(M)} x$ (compare with
Proposition 2.1).

An operator $x\in S(M)$ is said to be $\Phi$-{\it integrable} if
there exists a sequence $\{x_{n}\}\subset M$ such that
$x_{n}\xrightarrow{\Phi} x$ and $\Phi (\mid x_{n}
-x_{m}\mid)\xrightarrow{t(B)}0$ at $n, m \rightarrow \infty. $ It
follows from inequalities $\mid\Phi (x_{n}) - \Phi (x_{m})\mid
\leq\Phi (\mid x_{n} -x_{m}\mid)$ and completeness of the
topological $^*$-algebra $(S(M), t(M))$ that there exists an
element $\hat{\Phi}(x)\in S(B)$ such that $\Phi (x_{n})
\xrightarrow{t(B)}\hat{\Phi}(x).$ It is shown in \cite{CHZ} that
this limit $\hat{\Phi}(x)$ does not depend on the choice of a
sequence $\{x_{n}\}\subset M,$ for which $x_{n}\xrightarrow{\Phi}
x$ and $\Phi (\mid x_{n} -x_{m}\mid)\xrightarrow{t(B)} 0.$ It is
clear, any operator $x$ from $M$ is $\Phi$-integrable and
$\hat{\Phi}(x) = \Phi (x).$

Denote by $L^{1}(M, \Phi)$ the set of all $\Phi$-integrable
operators from $S(M),$ and for each $x\in L^{1}(M, \Phi),$ set
$\parallel x\parallel_{\mathbf{1}, \Phi} = \hat{\Phi}(\mid
x\mid).$ It is proved in [15] that $L^{1}(M, \Phi)$ is a linear
space in $S(M),$ in addition the following statement holds:\\

\noindent\textbf{Theorem 2.4} (\cite{ZCH, CHZ}).

$(i)$ {\it The mapping $\hat{\Phi}: L^{1}(M, \Phi)\rightarrow
S(B)$ has the following properties:

1) $\hat{\Phi}$ is a linear positive mapping, in particular,
$\hat{\Phi}(x^*) = (\hat{\Phi})^{*};$

2) $\hat{\Phi}(xy) = \hat{\Phi}(yx)$ and $\mid\hat{\Phi}(xy)\mid
\leq \parallel x\parallel_{M}\hat{\Phi}(y)$ for any $x\in M,$
$y\in L^{1}(M, \Phi);$

3) $x\in L_{1}(M,\Phi)\Leftrightarrow |x|\in L_{1}(M,\Phi),$
moreover $ \hat{\Phi}(\mid x^{*}\mid) = \hat{\Phi}(\mid x\mid)$ and
$\hat{\Phi}(\mid x\mid) = 0 \ \Leftrightarrow \ x=0;$

4) $\hat{\Phi}(\mid x+y\mid) \leq \hat{\Phi}(\mid x\mid) +
\hat{\Phi}(\mid y\mid);$

5) If  $\mid y\mid \leq \mid x\mid,$ $y\in S(M),$ $x\in L^{1}(M,
\Phi),$ then $y\in L^{1}(M, \Phi)$ and $\|y\|_{1,\Phi}\leq
\|x\|_{1,\Phi}$.

$(ii)$ $S(A)\cdot L^{1}(M, \Phi)\subset L^{1}(M, \Phi),$ moreover,
$\hat{\Phi}(zx) = \psi (z)\hat{\Phi}(x)$ for all $z\in S(A),$
$x\in L^{1}(M, \Phi)$, where $\psi$ is an extension of the
$\ast$-isomorphism from Theorem 2.3 to an $\ast$-isomorphism from
$S(A)$ onto $S(B)$.

$(iii)$ $(L^{1}(M, \Phi), \parallel\cdot\parallel_{\mathbf{1},
\Phi})$ is a Banach-Kantorovich space.}\\

For any $p > 1,$ set $L^{p}(M, \Phi) = \{x\in S(M): \ \mid
x\mid^{p}\in L^{1}(M, \Phi)\}$ and $\parallel x\parallel_{p, \Phi}
=\hat{\Phi}(\mid x\mid^{p})^{1/p}$ if $ x\in L^{p}(M, \Phi).$

We need the following properties of the spaces $(L^{p}(M, \Phi),
\parallel x\parallel_{p, \Phi})$ from \cite{CZ1}.\\

\noindent \textbf{Theorem 2.5} (\cite{CZ1}) $(i)$ {\it If an
element $x$ belongs to $L^{p}(M, \Phi),$ then the elements $x^{*}$
and $\mid x\mid$ also belong to $L^{p}(M, \Phi)$ and $\parallel
x^{*}\parallel_{p, \Phi} = \parallel \mid x\mid\parallel_{p, \Phi}
= \parallel x\parallel_{p, \Phi};$

$(ii)$ If $p, q >1,$ $\frac{1}{p} + \frac{1}{q} =1,$ $x\in
L^{p}(M, \Phi),$ $y\in L^{q}(M, \Phi)$ then $xy\in L^{1}(M, \Phi)$
and $\parallel xy\parallel_{1, \Phi} \leq \parallel x\parallel_{p,
\Phi} \parallel y\parallel_{q, \Phi};$

$(iii)$ $ML^{p}(M, \Phi) M \subset L^{p}(M, \Phi)$ and $\parallel
axb\parallel_{p, \Phi} \leq \parallel a\parallel_{M} \parallel
b\parallel_{M} \parallel x\parallel_{p, \Phi}$ for all $a, b\in M,
\ x\in L^{p}(M, \Phi);$

$(iv)$ If $\mid y\mid \leq\mid x\mid,$ $y\in S(M),$ $x\in L^{p}(M,
\Phi),$ then  $y\in L^{p}(M, \Phi)$ and $\parallel y\parallel_{p,
\Phi} \leq \parallel x\parallel_{p, \Phi};$

$(v)$ $L^{p}(M, \Phi)$ is a linear subspace in $S(M),$ moreover, $
M\subset L^{p}(M, \Phi)$ and $(L^{p}(M, \Phi),
\parallel\cdot\parallel_{p, \Phi})$ is a Banach-Kantorovich space;

$(vi)$ $S(A)L^{p}(M,\Phi)\subset L^{p}(M,\Phi),$ moreover,
$\parallel zx\parallel_{p,\Phi}=\psi(z)\parallel
x\parallel_{p,\Phi}$ for all $z\in S(A),$ $x\in L^{p}(M,\Phi),$
where $\psi$ is the $\ast$-isomorphism from Theorem 2.4} $(ii);$

$(vii)$ If  $\{x_{\alpha}\}\subset  \ L^{p}_{+}(M, \Phi),$
$\{x_{\alpha}\} \downarrow 0,$ then $\parallel x\parallel_{p,
\Phi} \downarrow 0.$

\section{The Arens algebras $\mathbf{L^{\omega}(M,\Phi)}$}

Let $\Phi$ be a Maharam trace on a von Neumann algebra $M$ with
values in $S(B)$. To define the Arens algebras associated with the
trace $\Phi,$ we need the following version of the H\"{o}lder
inequality.\\

\noindent \textbf{Theorem 3.1.} {\it If $ p,q,r>1,$ $\frac{1}{p} +
\frac{1}{q} = \frac{1}{r},$ $x \in L^{p}(M,\Phi),$ $y \in
L^{q}(M,\Phi)$, then $xy \in L^{r}(M,\Phi)$ and} $\parallel xy
\parallel _{r,\Phi} \leq \parallel x \parallel _{p,\Phi} \parallel y \parallel_{q,\Phi}.$

\begin{proof} Let $Q$ be the Stone compact corresponding to a
complete Boolean algebra $P(B)$ of all project from $B$. Identify
the algebra $S_{h}(B)$ with the algebra $C_{\infty}(Q)$ of all
continuous functions $f:Q\rightarrow [-\infty, +\infty]$ taking
values $\pm\infty$ only on nowhere dense sets from $Q$. As well as
in \cite{CZ1}, one can show that the element $\Phi( \mathbf{1})$
is reversible in $S(B),$ and a finite trace is defined on $M$ for
each $t\in Q$ by the equality $\varphi_{t}(x)=(\Phi(
\mathbf{1})^{-1}\Phi(x))(t)$. According to (\cite{D12}, 6.2.2),
the set $N_{t}=\{x\in M:\varphi_{t}(x^{\ast}x)=0\}$ is a two-sided
$\ast$-ideal in $M$. Consider the factor space $M/N_{t}$ with the
scalar product $([x],[y])=\varphi_{t}(y^{\ast}x)$ where $[x],[y]$
are the equivalence classes from $M/N_{t}$ with representatives
$x$ and $y,$ respectively.

Denote by $(H_{t},(\cdot,\cdot)_{t})$ the Hilbert space being the
completion of $(M/N_{t},(\cdot,\cdot)_{t}).$ Define the
$\ast$-homomorphism $ \pi_{t} :M \rightarrow B(H_{t}),$ setting
$\pi_{t}(x)([y])=[xy],$ $x,y\in M$. Denote by $U_{t}(M)$ the von
Neumann algebra in $B(H_{t})$ generated by the operators
$\pi_{t}(x),$ $x\in M$.

According to (\cite{D11}, 6.2), there exists a faithful normal
semi-finite trace $\tau_{t}$ on $(U_{t}(M))_{+}$ such that
$$
\tau_{t}(\pi_{t}(x^{2}))=([x],[x])=\varphi_{t}(x^{\ast}x)
$$
for all $x\in M_{+}$. Hence, $\tau_{t}(\pi_{t}(y))=\varphi_{t}(y)$
for all $y\in M_{+},$ moreover, $\tau_{t}(\mathbf{1}_{B(H_{t})})
=\varphi_{t}(\mathbf{1})<\infty,$ i.e. $\tau_{t}$ is a faithful
normal finite trace on $U_{t}(M).$

Let $L^{p}(U_{t}(M),\tau_{t})$ be a noncommutative $L^{p}$-space
associated with a numerical trace $\tau_{t}.$ According to
(\cite{FK1}, Theorem 4.9), we have
$$
\tau_{t}(|\pi_{t}(xy)|^{r})^{\frac{1}{r}} \leq\tau_{t}(|\pi_{t}
(x)|^{p})^{\frac{1}{p}}\tau(|\pi_{t}(y)|^{q})^{\frac{1}{q}}
$$
for all $x,y \in M.$

Since $\pi_{t}(\mid x\mid ) = \mid\pi_{t}(x)\mid $ for all $x \in
M,$ then by virtue of (\cite{D11},1.5.3) we have $\pi_{t}(\mid
x\mid^{p})=(\pi_{t}(\mid x\mid))^{p}.$ Hence,
$$
\varphi_{t}(\mid xy\mid^{r})^{\frac{1}{r}}\leq\varphi_{t}( \mid
x\mid^{p})^{\frac{1}{p}}\varphi_{t}(\mid y\mid^{q})^{\frac{1}{q}}
$$
or
$$[(\Phi(\mathbf{1}))^{-1}\Phi(\mid
xy\mid^{r}))(t)]^{\frac{1}{r}}\leq
[((\Phi(\mathbf{1}))^{-1}\Phi(\mid
x\mid)^{p})(t)]^{\frac{1}{p}}[((\Phi(\mathbf{1}))^{-1}\Phi(\mid
y\mid)^{q}))(t)]^{\frac{1}{q}}
$$
for all $t\in Q$ and $x,y\in M.$ It means that
$$
((\Phi(\mathbf{1}))^{-1}\Phi(\mid xy\mid^{r}))^{\frac{1}{r}}\leq
((\Phi(\mathbf{1}))^{-1}\Phi(\mid
x\mid^{p}))^{\frac{1}{p}}((\Phi(\mathbf{1}))^{-1}\Phi(\mid
y\mid^{q}))^{\frac{1}{q}}.
$$

Multiplying the both parts of the last inequality by
$\Phi(\mathbf{1}),$ we obtain
$$
\parallel xy \parallel_{r,\Phi}\leq\parallel x \parallel_{p,\Phi}\parallel y
\parallel_{q,\Phi}
$$
for any $x,y\in M.$

Let now $x\in L^{p}(M,\Phi),$ $y\in L^{q}(M,\Phi).$ Let us show
that $xy\in L^{r}(M,\Phi)$ and $\parallel xy
\parallel_{r,\Phi}\leq \parallel x \parallel_{p,\Phi}\parallel y
\parallel_{q,\Phi}.$ Set $a_{n}=E_{n}(|x|)|x|, b_{n}=E_{n}(|y|)|y|.$
We have $a_{n},b_{n}\in M_{+}$ and $a_{n}\uparrow|x|,
b_{n}\uparrow|y|,$ moreover,  $a_{n}\xrightarrow{\Phi}|x|,$
$b_{n}\xrightarrow{\Phi}|y|.$

Let $x=u\mid x\mid$ ($y=v\mid y \mid$, respectively) be the polar
decomposition for $x$ (for $y$, respectively) with the unitary
$u\in M ( v\in M).$ It is clear, $ua_{n}\xrightarrow{\Phi}x,$
$vb_{n}\xrightarrow{\Phi}y,$ and therefore
$ua_{n}vb_{n}\xrightarrow{\Phi}xy,$ in addition, $ua_{n}vb_{n}\in
M$ for all $n.$ Since the operations $z \mapsto |z|,$ $v\mapsto
v^{r},$ $v\geq 0$ are continuous in the topology $t(M)$, we get
$|ua_{n}vb_{n}|^{r}\xrightarrow{t(M)}|xy|^{r}$. The Fatou theorem
[11, Theorem 3.2 (iv)] implies that $xy\in L^{\omega}(M,\Phi)$ and
$$
\parallel xy \parallel^{r}_{r,\Phi}= \widehat{\Phi}(\mid
xy\mid^{r})\leq \sup\limits_{n\geq 1}\widehat{\Phi}(\mid
ua_{n}vb_{n}\mid^{r}) = \sup\limits_{n\geq 1}
\{\parallel(ua_{n})(vb_{n})\parallel^{r}_{r,\Phi}\}\leq
$$
$$
\leq\sup\limits_{n\geq 1} \parallel
ua_{n}\parallel^{r}_{p,\Phi}\parallel
vb_{n}\parallel^{r}_{q,\Phi}= \sup\limits_{n\geq1} \parallel
a_{n}\parallel_{p,\Phi}^{r}\parallel b_{n}\parallel_{q,\Phi}^{r}
\leq \parallel x\parallel^{r}_{p,\Phi}\parallel y\parallel^{r}_{q,
\Phi},
$$
i.e. $\parallel xy \parallel_{r, \Phi} \leq
\parallel x\parallel_{p,\Phi}\parallel
y\parallel_{q,\Phi}.$ \end{proof}

Suppose $L^{\omega}(M, \Phi) = \bigcap\limits_{p\geq 1} L^{p}(M,
\Phi).$ Theorems 2.5 and 3.1 imply the following\\

\noindent \textbf{Corollary 3.2.} $(i)$ \ $L^{\omega}(M,\Phi)$
{\it is a $\ast$-subalgebra in $S(M),$ and $M\subset
L^{\omega}(M,\Phi);$

\noindent $(ii)$ If $y\in S(M),$ $x\in L^{\omega}(M,\Phi),$
$|y|\leq |x|$, then $y\in L^{\omega}(M,\Phi)$;

\noindent $(iii)$  $S(A)\subset S(A)L^{\omega}(M,\Phi)\subset
L^{\omega}(M,\Phi)$ where $A$ is a $\ast$-subalgebra in $Z(M)$,
$\ast$-isomorphic to $B$ (see Theorem 2.3), i.e.
$L^{\omega}(M,\Phi)$ is left and right $S(A)$-module.} \\

An $\ast$-algebra $L^{\omega}(M,\Phi)$ is said to be \textit{the
Arens algebra} associated with the von Neumann algebra $M$ and the
Maharam trace $\Phi.$ If $\dim M <\infty$, then $M=S(M),$ and
therefore $M=L^{\omega}(M,\Phi).$ The converse assertion is also true.\\

\noindent \textbf{Proposition 3.3.} {\it If
$M=L^{\omega}(M,\Phi)$, then} $\dim M <\infty.$

\begin{proof} Suppose that $\dim M=\infty$. Then there exists a
countable set of mutually orthogonal nonzero projections
$\{q_{n}\}_{n=1}^{\infty} \subset P(M)$ for which
$\sup\{q_{n}\}=\mathbf{1}.$ Since $\Phi(\mathbf{1})=
\sum\limits_{n=1}^{\infty}\Phi(q_{n})$, the sequence $\Phi(q_{n})$
$(o)$-converges to zero in $S_{h}(B).$ Consider the spectral
family of projections $\{E_{\lambda}(\Phi(q_{n}))\}_{\lambda>0}$
and, using the inequality $\lambda E_{\lambda}^{\perp}
(\Phi(q_{n}))\leq \Phi(q_{n}),$ choose the sequence of numbers
$n_{1}< \dots < n_{k}< \dots$ such that $e_{k}=E_{\frac{1}{k^{2}}}
(\Phi(q_{n_{k}}))\neq 0, k=1,2\dots$ Let $\psi:A\rightarrow B$ be
the $^{*}$-isomorphism from Theorem 2.3. Set $x=\sum\limits_{k=
1}^{\infty}(\ln(k))q_{n_{k}}\psi^{-1}(e_{k})$ (the series
converges in the topology $t(M)$). It is clear, $x\in S_{+}(M)$
and $\Phi(x^{p})=\sum\limits_{k=1}^{\infty}
(\ln(k)^{p})\Phi(q_{n_{k}})e_{k} \leq\sum\limits_{k=1}^{\infty}
\frac{(\ln(k))^{p}}{k^{2}}e_{k}\leq (\sum\limits_{k=1}
^{\infty}\frac{(\ln(k))^{p}}{k^{2}})\mathbf{1_{B}}\in B,$ i.e.
$x\in L^{p}(M,\Phi)$ for all $p\geq 1.$ It means that $x\in
L^{\omega}(M,\Phi).$ Since $\Phi(q_{n_{k}}\psi^{-1}(e_{k}))=
e_{k}\Phi(q_{n_{k}})\neq 0,$ we get $q_{n_{k}}\psi^{-1}(e_{k})\neq
0$ for all $k=1,2,\dots,$ i.e. $x \notin M.$
\end{proof}

\noindent \textbf{Remark 3.4.} {\it If $M=Z(M)=A$, then by
Corollary 3.2}, $(iii),$ $S(M)=L^{\omega}(M,\Phi)$.\\

Denote by  $\mathbf{U}$ the base of zero neighborhoods in
$(S(B),t(B))$ consisting of ideal sets in the form (1). For any
$V\in \mathbf{U},$ $p\geq 1,$ set
$$
W(V,p)=\{x\in L^{p}(M,\Phi): \ \parallel x\parallel_{p,\Phi}\in
W\}.
$$

According to [24, chapter $I,\S 1$], there exists a topology
$\tau_{p}(M,\Phi)$  in $L^{p}(M,\Phi),$ with respect to which
$L^{p}(M,\Phi)$ is a Hausdorff topological vector space, and the
system of sets
$$\{x + W(V,p): V\in \mathbf{U}\}$$
forms the base of neighborhoods of the operator $x\in
L^{p}(M,\Phi).$ Convergence of the net $\{x_{\alpha}\}\subset
L^{p}(M,\Phi)$ to the operator $x\in L^{p}(M,\Phi)$ (notation:
$x_{\alpha}\xrightarrow{\tau_{p}(M,\Phi)} x,$) means that $\|
x_{\alpha}- x\|_{p,\Phi}\xrightarrow{t(B)} 0$.

Since $L^{p}(M,\Phi)$ is a Banach-Kantorovich space (Theorem 2.5
$(v)$), we have that $L^{p}(M,\Phi)- \tau_{p}(M,\Phi)$ is complete
\cite{CHY}, i.e. any $\tau_{p}(M,\Phi)$-fundamental net from
$L^{p}(M,\Phi)$ converges in $(L^{p}(M,\Phi),\tau_{p}(M,\Phi))$.

Now consider the set $W_{\omega}(V,p)=L^{\omega}(M,\Phi)\cap
W(V,p)$ and denote by $\tau_{\omega}(M,\Phi)$ the  Hausdorff vector
topology in $L^{\omega}(M,\Phi),$ in which the system of sets
$$\{x+W(V,p):V\in\mathbf{U}, p\geq 1\}$$
forms the base of neighborhoods of the element $x\in
L^{\omega}(M,\Phi).$ For a net$ \{x_{\alpha}\}\in
L^{\omega}(M,\Phi),$ its convergence $x_{\alpha}
\xrightarrow{\tau_{\omega}(M,\Phi)} x \in L^{\omega}(M,\Phi)$
means that $\|x_{\alpha}-x\|_{p,\Phi}\xrightarrow{t(B)} 0$ for all
$p\geq 1.$

If $1\leq r< p< \infty,$ $q = \frac{rp}{p-r},$ then
$\frac{1}{p}+\frac{1}{q}=\frac{1}{r},$ and by Theorem 3.1, we have
\begin{equation}
\|x\|_{r,\Phi}\leq
\|x\|_{p,\Phi}\|\mathbf{1}\|_{q,\Phi}=\Phi(\mathbf{1})^{\frac{p-r}{rp}}\|x\|_{p,\Phi}
\end{equation}
for all $x\in L^{\omega}(M,\Phi).$  It means that the topology
$\tau_{\omega}(M,\Phi)$ has the base of zero neighborhoods
consisting of sets in the form of $W_{\omega}(V,n)$, where $ V\in
\mathbf{U},$ $n\in \mathbb N,$ $\mathbb N$ is the set of all
natural numbers, i.e.
the following statement is valid:\\

\noindent \textbf{Proposition 3.5.} {\it If $B$ is a
$\sigma$-finite von Neumann algebra, then $\tau_{\omega}(M,\Phi)$
is a metrizable topology.} \\

Let $A$ be a von Neumann subalgebra in $Z(M)$ and let $\psi$ be an
$\ast$-isomorphism from $S(A)$ onto $S(B)$ being the extension of
the isomorphism from Theorem 2.3. If $x_{\alpha}, x\in S(A)$, then
$x_{\alpha},x\in L^{\omega}(M,\Phi)$ (Corollary 3.2 $(iii)$), and
by virtue of the equality $\|x_{\alpha}-x\|_{P,\Phi}=
\psi(x_{\alpha}-x)\|\mathbf{1}\|_{P,\Phi}$ (Theorem 2.5 $(vi)$),
the convergence $x_{\alpha}\xrightarrow{\tau_{\omega}(M,\Phi)} x$
is equivalent to the convergence $\psi(x_{\alpha})
\xrightarrow{t(B)} \psi(x)$. Hence,
\begin{equation}
x_{\alpha}\xrightarrow{\tau_{\omega}(M,\Phi)} x
\Longleftrightarrow x_{\alpha}\xrightarrow{t(A)} x,
\end{equation}
i.e. the topology $\tau_{\omega}(M,\Phi)$ induces on $S(A)$ the
topology $t(A)$. Therefore metrizability of the topology
$\tau_{\omega}(M, \Phi)$ implies metrizability of the topology
$t(A),$ that is equivalent to $\sigma$-finiteness of the von
Neumann algebra $B.$ Thus, taking into account Proposition 3.5, we
obtain the following criterion for metrizability of
$\tau_{\omega}(M,\Phi)$.\\

\noindent \textbf{Theorem 3.6.} {\it The following conditions are
equivalent for the Arens algebra $L^{\omega}(M,\Phi)$ associated
with the von Neumann algebra $M$ and an $S(M)$-valued trace
$\Phi:$

$(i)$. The topology $\tau_{\omega}(M,\Phi)$ is metrizable;

$(ii)$. The von Neumann algebra $B$ is $\sigma$-finite;

$(iii)$. The von Neumann algebra $M$ is $\sigma$-finite}.

\begin{proof} Implications $(i)\Leftrightarrow (ii)$ are already proved.
Implication $(iii)\Longrightarrow(ii)$ is evident, since the
$\sigma$-finiteness of the algebra $B$ is equivalent to
$\sigma$-finiteness of the von Neumann algebra $A\subset Z(M)$ .
Implication $(ii)\Longrightarrow(iii)$ is proved as well as in
[1], where $\sigma$-finiteness of the algebra $M$ is obtained
using a center-valued trace $\Phi_{M}:M\longrightarrow Z(M)$ and
$\sigma$-finiteness of the center of $Z(M).$
\end{proof}

\noindent \textbf{Proposition 3.7.} $(L^{\omega}(M,\Phi),
\tau_{\omega}(M,\Phi))$ {\it is a complete topological
$\ast$-algebra.}

\begin{proof} Since  $(L^{p}(M,\Phi), \tau_{p}(M, \Phi))$ is complete
\cite{CHY}, for any $\tau_{\omega}(M,\Phi)$-fun\-da\-men\-tal net
$\{x_{\alpha}\}\subset L^{\omega}(M,\Phi)$ there exists $x_{p}\in
L^{\omega}(M,\Phi)$ such that $x_{\alpha}\xrightarrow{\tau_{p}
(M,\Phi)}x_{p}$. The inequality (3) implies that for $1\leq
q<p<\infty $ the topology $\tau_{q}(M,\Phi)$ is majorized on
$L^{\omega}(M,\Phi)$ by the topology $\tau_{p}(M,\Phi)$. Therefore
the convergence $x_{\alpha}\xrightarrow{\tau_{p}(M,\Phi)}x_{p}$
implies the convergence $x_{\alpha}\xrightarrow{\tau_{q} (M,\Phi)
}x_{p}$, what implies the equality $x_{p}=x_{q}$ at $1\leq
q<p<\infty $. Hence, $x:=x_{p}$ is an element from
$L^{\omega}(M,\Phi),$ and $x_{\alpha}\xrightarrow{\tau_{\omega}
(M,\Phi)}x_{p}.$ It means that $(L^{\omega}(M,\Phi),\tau_{\omega}
(M,\Phi))$ is complete. Further, the equality
$\|x\|_{p,\Phi}=\|x^{\ast}\|_{p,\Phi}$ (Theorem 2.5 $(i))$ and
inequality $\|xy\|_{p}\leq \|x\|_{2p}\|y\|_{2p}$ (Theorem 3.1),
imply that the involution operation is continuous in
$(L^{\omega}(M,\Phi),$ $\tau_{p}(M,\Phi)),$ and the multiplication
operation is continuous in both variables. Thus,
$(L^{\omega}(M,\Phi),\tau_{\omega}(M,\Phi))$ is a complete
topological $\ast$-algebra.
\end{proof}

Theorem 3.6 and Proposition 3.7 imply the following \\

\noindent \textbf{Corollary 3.8.} {\it If $B$ is a $\sigma$-finite
von Neumann algebra, then $(L^{\omega}(M,\Phi),
\tau_{\omega}(M,\Phi))$ is a complete metrizable topological
$\ast$-algebra, in particular, $(L^{\omega}(M,\Phi),
\tau_{\omega}(M,\Phi))$ is a $F$-space.}\\

Denote by $t_{\omega}(M)$ the topology on $L^{\omega}(M,\Phi)$
induced by the topology $t(M)$ from $S(M).$\\

\noindent \textbf{Proposition 3.9.} $(i)$ $t_{\omega}(M)\leq
\tau_{\omega}(M,\Phi);$

\noindent $(ii)$ {\it If $t_{\omega}(M)= \tau_{\omega}(M,\Phi)$,
then} $L^{\omega}(M,\Phi)=S(M).$

\begin{proof} $(i).$ Let $\{x_{\alpha}\}\subset
L^{\omega}(M,\Phi),$ and $x_{\alpha}\xrightarrow{\tau_{\omega}
(M,\Phi)}0,$ i.e. $\|x_{\alpha}\|_{p,\Phi}\xrightarrow{t(B)}0$ for
all $p\geq 1.$

Since $ \parallel|x|\parallel = \parallel x\parallel_{p,\Phi}$
(Theorem 2.5 $(i)$), we obtain $\parallel|x_{\alpha}|
\parallel_{p,\Phi} \xrightarrow{t(B)}0.$ Let $\{E_{\lambda}(|
x_{\alpha}|)\}_{\lambda > 0}$ be the spectral family of
projections for $\mid x_{\alpha}\mid.$ By virtue of the inequality
$$
\Phi(E_{\lambda}^{\perp}(\mid x_{\alpha}\mid))\leq
\frac{1}{\lambda}\||x|_{\alpha}\|_{\mathbf{1},\Phi},\ \lambda
>0,
$$
we have $x_{\alpha}\xrightarrow{\Phi}0,$ and therefore
$x_{\alpha}\xrightarrow{t(M)}0 \ \cite{ZCH}$.

$(ii).$ Suppose that $t_{\omega}(M)=\tau_{\omega}(M,\Phi).$ Since
$(L^{\omega}(M,\Phi), \tau_{\omega}(M,\Phi))$ is complete,
$L^{\omega}(M,\Phi)$ is a closed subalgebra in $(S(M),t(M)).$ Let
$x\in S_{+}(M),$ $x_{n}=E_{n}(x)x,$ $n\in N.$ It is clear that
$x_{n}\in M,$ and $\lambda E^{\perp}_{\lambda}(x-x_{n})\leq
(x-x_{n})\downarrow 0,$ that's why $ \Phi(E^{\perp}_{\lambda}
(x-x_{n}))\xrightarrow{t(B)}0$ at $n\rightarrow\infty$ for all
$\lambda>0$. The means that $x_{n}\xrightarrow{\Phi}x $, and
therefore $x_{n}\xrightarrow{t(M)} x$. Hence, $x\in
L^{\omega}(M,\Phi)$. Since each element from $S(M)$ is a linear
combination of four elements from $S_{+}(M)$, we have $x\in
L^{\omega}(M,\Phi)$ for all $x\in S(M),$ which implies the
equality $L^{\omega}(M,\Phi)=S(M)$. \end{proof}

\noindent \textbf{Remark 3.10.} {\it If $\dim M<\infty,$ then
$L^{\omega}(M,\Phi)=M=S(M)$, and the both vector topologies
$\tau_{\omega} (M,t)$ and $t(M)$ coincide on $L^{\omega}(M,\Phi),$
moreover they are equal to the topology generated by the}
$C^{*}$-norm $\|\cdot\|_{M}.$

\section{Disjunct completeness of Arens algebras}

Let $L^{\omega}(M,\Phi)$ be an Arens algebra associated with the
von Neumann algebra $M$ and $S(B)$-valued Maharam trace $\Phi$,
and let  $A$ be a von Neumann subalgebra in $Z(M),$ and $\psi$ be
an $\ast$-isomorphism from $S(A)$ onto $S(B)$ from Theorem 2.4
$(ii)$. Consider arbitrary decomposition $\{e_{i}\}_{i\in I}$ of
the unit $\mathbf{1}_{B}$ of the complete Boolean algebra $P(B)$
of all projections from the commutative von Neumann algebra $B.$
It is clear that $q_{i}=\psi^{-1}(e_{i}),$ $i\in I$, is
decomposition of the unit $\mathbf{1}_{A}=\mathbf{1}_{M}$ of the
Boolean algebra $P(M)$. For each $q_{i}x\in q_{i}M,$ we have
$\Phi(q_{i}x)= \psi(q_{i}) \Phi(x)=e_{i}\Phi(x)\in
e_{i}S(B)=S(e_{i}B)$. Hence, the restriction $\Phi_{i}$ of the
Maharam trace $\Phi$ on $e_{i}M$ is an $S(e_{i}B)$-valued Maharam
trace on the von Neumann algebra $e_{i}M$ for all $i\in I$. If
$p\geq1$ and $x\in L^{p}(M,\Phi)$, then, evidently, $q_{i}x\in
L^{p}(q_{i}M,\Phi_{i}),$ and
$e_{i}\|x\|_{P,\Phi}=\|q_{i}x\|_{P,\Phi_{i}}$ for all $i\in I$.
Since $(L^{p}(M,\Phi),\|\cdot\|_{P,\Phi})$ is a Banach-Kantorovich
space (Theorem 2.5), then $L^{p}(M,\Phi)$ is disjunct complete
(\cite{K11}, 2.1.5, 2.2.1), i.e. the following assertion is valid.\\

\noindent \textbf{Proposition 4.1.} {\it If $x_{i}\in
L^{p}(q_{i}M,\Phi_{i})$ for all $i\in I$, a there exists the
unique element $x\in L^{p}(M,\Phi)$, for which $q_{i}x=x_{i}$ at
all} $i\in I$. \\

The following property of disjunct completeness for the Arens
algebra $L^{\omega}(M,\Phi)$ follows immediately from Proposition 4.1.\\

\noindent \textbf{Corollary 4.2} {\it Let $\{e_{i}\}_{i\in I}$ be
arbitrary decomposition of the unit of the Boolean algebra $P(B)$,
$q_{i}=\psi^{-1}(e_{i})$, $x_{i}\in L^{\omega}(q_{i}M,\Phi_{i})$,
$i\in I$. Then there exists a unique $x\in L^{\omega}(M,\Phi)$
such that $q_{i}x=x_{i}$ for all} $i\in I$.\\

Consider the direct product $\prod\limits_{i\in I}
L^{\omega}(q_{i}M,\Phi_{i})$ of $\ast$-algebras
$L^{\omega}(q_{i}M,\Phi_{i})$  with coordinate-wise algebraic
operations and involution. Define the mappings
$$
U: L^{\omega}(M,\Phi)\longrightarrow\prod\limits_{k\in
I}L^{\omega}(q_{i}M,\Phi_{i}),
$$
setting $U(x)=\{q_{i}x\}_{i\in I} $. According to Corollary 4.2,
the mapping $U$ is an $\ast$-isomorphism from $L^{\omega}(M,\Phi)$
onto the $\ast$-algebra $\prod\limits_{i\in I}L^{\omega}
(q_{i}M,\Phi_{i})$. Denote by $t\left(\{q_{i}\} \right)$ the
Tychonoff product of topologies $\tau_{\omega}(q_{i}M,\Phi_{i})$
in $\prod\limits_{i\in I}L^{\omega}(q_{i}M,\Phi_{i})$. By virtue
of properties of Tychonoff topologies, the pair $\left(
\prod\limits_{i \in I} L^{\omega} \left(q_{i}M,\Phi_{i}\right),
t\left(\{q_{i}\}\right)\right)$ is a   complete topological $\ast$-algebra.\\

\noindent \textbf{Proposition 4.3.} {\it The mapping
$$
U:(L^{\omega}(M,\Phi),\tau_{\omega}(M,\Phi))\longrightarrow
\left(\prod\limits_{i\in I}L^{\omega}
(q_{i}M,\Phi_{i}),t(\{q_{i}\})\right)
$$
is a homeomorphism, in particular,
$x_{\alpha}\xrightarrow{\tau_{\omega}(M,\Phi)}x, x_{\alpha},x\in
L^{\omega}(M,\Phi)$ iff $q_{i}x_{\alpha} \xrightarrow{\tau_{
\omega} (q_{i}M,\Phi_{i})} q_{i}x$ for all} $i\in I$.\\

The proof of this statement follows from the definition of the
Tychonoff topology and equalities $\|q_{i}x\|_{P,\Phi}=e_{i}
\|x\|_{P,\Phi_{i}}$ for all $i\in I, p\geq 1, x\in
L^{\omega}(M,\Phi)$ (see Theorem 2.5 $(vi)$).

The following theorem gives necessary and sufficient conditions
for locally convexity of the topology $\tau_{\omega}(M,\Phi).$ \\

\noindent \textbf{Theorem 4.4.} {\it Let $L^{\omega}(M,\Phi)$ be
an Arens algebra associated with a von Neumann algebra $M$ and an
$ S(B)$-valued Maharam trace $\Phi$. The following conditions are
equivalent:

$(i).$ The topology $\tau_{\omega}(M,\Phi)$ is locally convex;

$(ii).$ $B$ is an atomic von Neumann algebra}.

\begin{proof} $(i)\Rightarrow (ii)$. Let
$\tau_{\omega}(M,\Phi)$ is a locally convex topology. According to
(4), the topology $t(A)$ is also locally convex. Hence, by
Proposition 2 from  (\cite{SAH}, chapter $V,\S 3$], the von
Neumann algebra $A$ is atomic, what implies atomicity of the von
Neumann algebra $B$.

$(ii)\Rightarrow (i)$. If $B$ is an atomic von Neumann algebra,
then there exists in $P(B)$ a decomposition $\{e_{i}\}_{i\in I}$
of the unit such that $e_{i}$ is an atom in $P(B)$ for each $i\in
I$. Therefore $e_{i}S(B)=\mathbb{C}$ and $S(B)=\mathbb{C}^{I}$. In
this case, for $g_{i}=\psi^{-1}(e_{i}),$ the function $\Phi_{i}
(g_{i}x)=e_{i}\Phi(x)$ is a faithful normal finite trace on the
von Neumann algebra $g_{i}M,$ and thus the topology $\tau_{\omega}
(g_{i}M,\Phi_{i})$ is locally convex \cite{Z11}. Hence, the
Tychonoff topology is also locally convex. It remains to use
Proposition 3.4, by virtue of which, the topology
$\tau_{\omega}(M,\Phi)$ is also locally convex.
\end{proof}

\noindent\textbf{Remark 4.5.} {\it If the topology $\tau_{\omega}
(M,\Phi)$ is normable, then the topology $t(A)$ is also normable
(see (4)), what implies finite dimensionality of the algebra $A$
(see Proposition 4 from [27, chapter $V,\S 3$]). Hence, the
algebra $B$ is also finite-dimensional, and therefore there is a
finite collection of atoms $\{e_{i}\}_{i=1}^{n}$ â $P(B)$, for
which $\sum\limits_{i=1}^{n}e_{i}=\mathbf{1}_{B}.$ In this case
the $\ast$-algebra is $B-\ast$-isomorphic to $\mathbb{C}^{n},$ and
$\tau(x)=\sum\limits_{i=1}^{n}e_{i}\Phi(x)$ is a faithful normal
finite numeric trace on $M$, for which} $L^{\omega}(M,\Phi)=
L^{\omega}(M,\tau)$.

\section{Comparison of the $\mathbf{(o})$-topology and the to\-po\-lo\-gy $\mathbf{\tau_{\omega}(M,\Phi)}$}

Let $L^{\omega}(M,\Phi)$ be the Arens algebra associated with a
von Neumann algebra $M$ and an $S(B)$-valued Maharam trace $\Phi.$
Since the involution in $L^{\omega}(M,\Phi)$ is continuous with
respect to the topology $\tau_{\omega}(M,\Phi)$ (Proposition 3.7),
the set $L^{\omega}_{h}(M,\Phi)$ is closed in $(L^{\omega}
(M,\Phi),\tau_{\omega}(M,\Phi)).$ It was proved in \cite{Y12} that
the set $S_{+}(M)$ is closed in $(S(M),t(M))$. Therefore,
according to Proposition 3.9 $(i)$, the set
$L_{+}^{\omega}(M,\Phi)$ is closed in
$(L^{\omega}(M,\Phi),\tau_{\omega}(M,\Phi))$. Denote by
$\tau_{\omega h}(M,\Phi)$ the topology in $L^{\omega}_{h}(M,\Phi)$
induced by the topology $\tau_{\omega}(M,h)$ from
$L^{\omega}(M,\Phi),$ and by $\tau_{o}(M,\Phi)$ -- the
$(o)$-topology in $L^{\omega}_{h}(M,\Phi),$ i.e. the strongest
topology  in an ordered linear space $L^{\omega}_{h}(M,\Phi)$, for
which $(o)$-convergence of nets
implies their topological convergence.\\

\noindent \textbf{Theorem 5.1}.

$(i).$ $\tau_{\omega,h}(M,\Phi)\leq\tau_{o}(M,\Phi)$;

$(ii).$ $\tau_{\omega h}(M,h)=\tau_{o}(M,\Phi)$ {\it if and only
if a von Neumann algebra $B$ is $\sigma$-finite.}

\begin{proof} $(i).$  If $\{x_{\alpha}\}_{\alpha\in
A}\subset L^{\omega}_{h}(M,\Phi),$ and $\{x_{\alpha}\}$
$(o)$-converges in $L^{\omega}_{h}(M,\Phi)$ to an element $x\in
L^{\omega}_{h}(M,\Phi),$ then by definition there exist nets
$\{y_{\alpha}\}_{\alpha\in A},$ $\{z_{\alpha}\}_{\alpha\in A}$
from $L^{\omega}_{h}(M,\Phi)$ such that $y_{\alpha}\leq
x_{\alpha}\leq z_{\alpha}$ for all $\alpha\in A$ è
$y_{\alpha}\uparrow\ x, z_{\alpha}\uparrow\ x$. Since
$y_{\alpha}-x\leq x_{\alpha}-x\leq z_{\alpha}-x$, we have
$$
0\leq (x_{\alpha}-x)+(x-y_{\alpha})\leq (z_{\alpha}-x)+(x-y_{\alpha})
$$

Using convergences $(z_{\alpha}-x)\downarrow 0,$
$(x-y_{\alpha})\downarrow 0,$ and Theorem 2.5 $(iv),$ $(vii)$, we
obtain
$$\|x_{\alpha}-x\|_{p,\Phi}\leq\|(x_{\alpha}-x)+(x-y_{\alpha})\|_{p,\Phi}+
\|x-y_{\alpha}\|_{p,\Phi}\leq
$$ $$
\le \|z_{\alpha}-x\|_{p,\Phi}+2\|x-y_{\alpha}\|_{p,\Phi}
\xrightarrow{t(B)}0
$$
for all $p\geq 1$. Therefore
$x_{\alpha}\xrightarrow{\tau_{\omega}(M,\Phi)}x$, which implies
of the inequality
$$\tau_{\omega,h}(M,\Phi)\leq \tau_{o}(M,\Phi).$$

$(ii).$\ If the von Neumann algebra $B$ is $\sigma$-finite, then
the topology $\tau_{\omega}(M,\Phi)$ is metrizable (Theorem 3.6).
Repeating the proof of Theorem 2 from \cite{A13}, we obtain
$\tau_{o}(M,\Phi)\leq \tau_{\omega h}(M,\Phi).$ Therefore,
according to $(i)$, the equality $\tau_{o}(M,\Phi)=\tau_{\omega
h}(M,\Phi)$ is valid.

Now suppose that $\tau_{o}(M,\Phi)=\tau_{\omega h}(M,\Phi)$. Let
$A$ and $\psi$ by the same as in Theorem 2.4. $(ii)$. Since
$S(A)\subset L^{\omega}(M,\Phi),$ and the topology
$\tau_{\omega,h}(M,\Phi)$ induces on $S(A)$ the topology $t(A)$
(see(4)), repeating the proof of Theorem 2 from \cite{CHM}, we
obtain that the $(o)$-topology in $S_{h}(A)$ coincides with the
topology $t_{h}(A)$, where $t_{h}(A)$ is the topology $S_{h}(M)$
induced by the topology $t(A)$ from $S(A)$. Hence, by Theorem 2
from \cite{CHM} the von Neumann algebra $A$ is $\sigma$-finite.
Since $\psi$ is $\ast$-isomorphism from $A$ onto $B$, the algebra
$B$ is also $\sigma$-finite.
\end{proof}

\end{document}